\documentclass[10pt,a4paper]{amsart}
\usepackage{amsmath,amssymb,amsthm}

\numberwithin{equation}{section}
\newtheorem{thm}{Theorem}[section]
\newtheorem{cor}[thm]{Corollary}
\newtheorem{prop}[thm]{Proposition}
\newtheorem{lemma}[thm]{Lemma}

\theoremstyle{remark}

\theoremstyle{definition}

\theoremstyle{remark}

\theoremstyle{remark}
\newtheorem{remark}{Remark}[section]

\newenvironment{Abstract}
{\begin{center}\textbf{\footnotesize{Abstract}}%
\end{center} \begin{quote}\begin{footnotesize}}
{\end{footnotesize}\end{quote}\bigskip}

\begin{document}

\newcommand\vrm{v_{R,M,T}}
\newcommand{\al}{\alpha}
\newcommand{\be}{\beta}
\newcommand{\de}{\partial}
\newcommand{\la}{\lambda}
\newcommand{\La}{\Lambda}
\newcommand{\ga}{\gamma}
\newcommand{\ep}{\epsilon}
\newcommand{\del}{\delta}
\newcommand{\Del}{\Delta}
\newcommand{\sig}{\sigma}
\newcommand{\ome}{\omega}
\newcommand{\Ome}{\Omega}
\newcommand{\uone}{u^{(1)}}
\newcommand{\utwo}{u^{(2)}}
\newcommand{\C}{{\mathbb C}}
\newcommand{\N}{{\mathbb N}}
\newcommand{\Z}{{\mathbb Z}}
\newcommand{\R}{{\mathbb R}}
\newcommand{\Rn}{{\mathbb R}^{n}}
\newcommand{\Rnu}{{\mathbb R}^{n+1}_{+}}
\newcommand{\Cn}{{\mathbb C}^{n}}
\newcommand{\spt}{\,\mathrm{supp}\,}
\newcommand{\Lin}{\mathcal{L}}
\newcommand{\SSS}{\mathcal{S}}
\newcommand{\F}{\mathcal{F}}
\newcommand{\xxi}{\langle\xi\rangle}
\newcommand{\xx}{\langle x\rangle}
\newcommand{\yy}{\langle y\rangle}
\newcommand{\dint}{\int\!\!\int}
\newcommand{\triple}[1]{{|\!|\!|#1|\!|\!|}}

\renewcommand{\Re}{\;\mathrm{Re}\;}


\title[The solution operator to the wave map system]%
{On the continuity of the solution operator to the wave map system}

\author{Piero D'Ancona}
\address{Piero D'Ancona: Unversit\`a di Roma ``La Sapienza'',
Dipartimento di Matematica, Piazzale A.~Moro 2, I-00185 Roma,
Italy}
\email{dancona@mat.uniroma1.it}

\author{Vladimir Georgiev}
\address{Vladimir Georgiev: Universit\`a di Pisa, Dipartimento di Matematica,
Via Buonarroti 2, I-56127 Pisa, Italy}
\email{georgiev@dm.unipi.it}

\thanks{The authors were partially supported
by the Programma Nazionale M.U.R.S.T.
``Problemi e Metodi nella Teoria delle Equazioni Iperboliche.''}
\keywords{AMS Subject Classification: %
35L70, 
58J45
. Keywords: Wave map, solution operator,
nonlinear equation, hyperbolic equation,
ill-posed problem}

\maketitle

\begin{Abstract}
We investigate the continuity properties of the solution operator to
the wave map system from
general nonflat target of arbitrary dimension, and we prove by an
explicit class of counterexamples that this map is not uniformly
continuous in the critical norms on any neighbourhood of 0.
\end{Abstract}

\section{Introduction}

Let $(N,g)$ be a smooth $d$-dimensional Riemannian manifold
with metric $g$, which with no
loss in generality we can
isometrically embed in $\R^m$ for some $m>d$.
For the functions $u:\R\times\R^{n}\to N$ defined on the flat
Minkowski space $\R_{t}\times\R^{n}_{x}$ with values in the target $N$,
consider the functional
$$J(u)=\int_{\R\times\R^{n}}
   \langle \de_{\al}u,\de^{\al}u\rangle_{g(u)} dtdx$$
where summation over $\al=0,1,...,n$ is intended, with
$$(\de_{0},...,\de_{n})=(\de_{t},\de_{x_{1}},...,\de_{x_{n}}),
   \qquad
   (\de^{0},...,\de^{n})=(\de^{t},-\de_{x_{1}},...,-\de_{x_{n}})$$
as usual, while $\langle\cdot,\cdot\rangle_{g}$ is the product in
the metric $g$.

The critical points of the functional $J$ are called \emph{wave maps}.
If we choose a system of coordinates on $N$, then locally smooth wave
maps satisfy the equation
\begin{equation}\label{wmsys}
         \square u^\ell + \Gamma^\ell_{bc}(u)
         \partial_\alpha u^b \partial^\alpha u^c = 0,\qquad
            \ell=1,\dots,m,
\end{equation}
where $\Gamma^{\ell}_{jk}$ denote the Christoffel
symbols on $N$ in the chosen coordinates. The natural problem for
this system of wave equations is clearly the Cauchy problem with data
at $t=0$
\begin{equation}\label{wmdata}
         u(0,x)=u_0,\qquad u_t(0,x)=u_1;
\end{equation}
the usual space for the data are Sobolev spaces
$$(u_{0},u_{1})\in H^{s}(\R^{n},N)\times H^{s-1}(\R^{n},N)$$
for suitable values of $s\in\R$. Here we used the space
\begin{equation}\label{Hs}
   H^s(\R^n;N)\equiv\{v\in H^s(\R^n;\R^m),
     \ v(\R^n)\subseteq N\},\qquad s\in\R
\end{equation}
with the induced norm; notice that $H^s(\R^n;N)=\emptyset$
if $0\not\in N$, but it causes no loss in generality to assume
that $0\in N$ after a translation in the ambient space.

An alternative description of the wave map system, which usually
gives a simpler expression in presence of symmetry of the target
is the following: a wave map is a function
$u:\R\times\R^{n}\to\R^{m}$ such that
$$ u(t,x)\in N, \qquad\square u\perp N\qquad\hbox{for all\ \ }(t,x).$$
A good introduction on this subject with comprehensive
references can be found in \cite{ShStr1}.

The Cauchy problem for wave maps has been extensively studied in recent years, starting with the work of Ginibre and Velo \cite{GinibrVelo82}.
Not many general results for the Cauchy problem are known.
Very roughly speaking and with several omissions,
one has the following types of results:
\begin{itemize}
\item Global existence of weak solutions when the target is compact
(several authors, see e.g. \cite{FMS}, \cite{Sch}, \cite{Yi}).
\item Local existence for data in $H^{s}$, $s>n/2$.
This is classical if $s$ is large enough, but for $s$ close to the
critical value $s=n/2$ it is a much more difficult result, due
to Klainerman, Machedon, Selberg and obtained through careful
bilinear estimates (see in particular \cite{KS}). See also
Tataru \cite{Tat} for the case of Besov spaces.
\item Global existence for small data. Again, this result can be proved
by ``standard'' methods in the smooth case (Y.~Choquet-Bruhat), but
the recent results of Tao (see in particular \cite{Ta2}, see also
\cite{ShStr2}) show that it is sufficient to assume that the data
are in $H^{s}$ for some $s>n/2$ and that they are small in the
homogeneous $\dot H^{n/2}$ norm.
\item In the presence of symmetry one has of course
sharper results;
to this class belong the radial case, considered by
Christodoulou, Tahvildar-Zadeh and more recently by Struwe
(see \cite{CT}, \cite{Str2}, \cite{Str3}), and the equivariant case,
for which a fairly complete theory exists, due to Shatah, Tahvildar-Zadeh, Struwe, and Grillakis (see e.g. \cite{G}, \cite{Sch},
\cite{ShT1}, \cite{Str4}).
\end{itemize}

The precise behaviour of the wave map system in the critical case
$s=n/2$ is an open problem (while in the \emph{subcritical} case
$s<n/2$ one has in general ill-posedness in a sharp sense, i.e.
non-uniqueness, see e.g. \cite{Ta1} and \cite{DG}). A possible line of
attack was suggested by Bourgain (see \cite{Bo1}, \cite{Bo2}
and also Tzvetkov \cite{Tz 1}) who proved
that the map data $\mapsto$ solution
to some nonlinear evolution equations is not $C^{2}$ in the subcritical
Sobolev spaces. This holds for the cubic NLS, for KdV and mKdV with
different critical indices. The result was sharpened by Kenig, Ponce and
Vega \cite{KPV} who proved that the solution map actually is
not (locally) uniformly
continuous in the subcritical spaces. We also mention \cite{BK}
and \cite{LST}
where the case of the supercritical nonlinear wave equation and
of the Benjamin-Ono equation are
considered.

Our aim here is to prove a similar result for the wave map system,
in the \emph{critical} case $s=n/2$.
Our assumption on $N$ will be quite general;
essentially we only require
that $N$ is not flat. More precisely, we assume that
\begin{equation}\label{assN}
   \hbox{\textit{there exists a geodesic curve $\ga:]-s_0,s_0[\to N$
   with $\ga(0)=0$,
   $\ga''(0)\neq0$.}}
\end{equation}
From such generality it will be clear that the ill-posedness in
the sense of uniform continuity is a general properties of
nonlinear equations like \eqref{wmsys} more than a geometric
property. Our result is the following:

\begin{thm}\label{thm1}
   Let $N$ be a smooth Riemannian manifold, isometrically
   embedded in $\R^m$,
   such that there esists a geodesic curve $\ga:]-s_0,s_0[\to N$ with
   \begin{equation}\label{gamma}
     \ga(0)=0\in N,\qquad \ga''(0)\neq0.
   \end{equation}
   Assume a solution map $\Phi:(u_0,u_1)\to u$ for system ~\eqref{wmsys}
   with data ~\eqref{wmdata}
   is defined on some neighbourhood $U$ of 0 in
   $X\times Y=H^{n/2}(\R^n;N)\times H^{n/2-1}(\R^n;N)$. Then, for
   any $T>0$,
   $\Phi$ is not uniformly continuous between the spaces
   $$\Phi:U\subseteq H^{n/2}\times H^{n/2-1}\to
       C([0,T];\dot H^{n/2}(\R^n,\R^m))$$
   or
   $$\Phi:U\subseteq H^{n/2}\times H^{n/2-1}\to
       C^1([0,T]; \dot H^{n/2-1}(\R^n,\R^m)).$$
\end{thm}

\goodbreak

As usual, we say that $\Phi:U\subseteq X\times Y\to L$
is uniformly continuous on $U$ if:
\textit{for any $\ep>0$ there exists $\del>0$ such that, for any
$(\uone_0,\uone_1)$ and $(\utwo,\utwo)$ in $U$}
\begin{equation}\label{unifcont}
   \|(\uone_0,\uone_1)-(\utwo_0,\utwo_1)\|_{X\times Y}\leq\del
   \quad\Rightarrow\quad
   \|\uone-\utwo\|_L\leq\ep
\end{equation}
\textit{where $\uone=\Phi(\uone_0,\uone_1)$, $\utwo=\Phi(\utwo_0,\utwo_1)$.}
Thus the above result excludes in particular that $\Phi$ is (locally)
Lipschitz or H\"older continuous.

\begin{remark}
It is not difficult to prove by similar arguments that the
solution map is not uniformly continuous also in the
\emph{subcritical} case, i.e., from $H^s\times H^{s-1}$, $1\leq
s<n/2$ with values in $C([0,T],H^s)$ or $C^1([0,T],H^{s-1})$.
However, it is already known, at least in the case of a
rotationally symmetric target, that a much stronger ill posedness
result holds, namely the local non-uniqueness can be proved. This
was obtained for $n=3$ in \cite{ShT1}, for $n\geq 4$ in
\cite{CST} and for $n=2$ in \cite{DG}. Since the arguments in
these results have a local nature, it is reasonable to argue that
non uniqueness may hold also in the general nonsymmetric case.
\end{remark}

\begin{remark}\label{remm1}
The proof of the Theorem is based on an explicit construction of
sequences of data such that the corresponding solutions violate
~\eqref{unifcont}; such solutions are of geodesic type, i.e., of the
form $\ga\circ v(t,x)$ with $v(t,x)$ a real valued solution of the
homogeneous wave equation. We recall
that if $\ga(s)=(\ga_{1},...,\ga_{m})$ is an arbitrary curve in $\R^{m}$ with values in $N$, and $v(t,x)$ an arbitrary real valued function, for the composition $u(t,x)=\ga(v(t,x))$ we can write
$$\square u^{\ell} + \Gamma^\ell_{bc}(u)
         \partial_\alpha u^b \partial^\alpha u^c \equiv
   \ga_{\ell}'\cdot\square v+
   \left(\ga''_{\ell}+ \Gamma^\ell_{bc}(\ga)\ga'_{b}\ga'_{c}\right)
      \cdot\partial_\alpha v \partial^\alpha v $$
and this is identically zero as soon as $\square v=0$ and $\ga(s)$
is a geodesic curve.
\end{remark}

\begin{remark}\label{rem0}
The ill posedness for the wave map problem in
the case $n=1$, $s=1/2$ is proved in \cite{Ta1}.
It is interesting to mention also the paper \cite{NaO}, where a \emph{scalar} wave equation of the form
$$\square u+f(u)\partial_{\al}u\partial^{\al}u=0$$
is studied in the critical spaces $H^{n/2}\times H^{n/2-1}$.
Of course in the scalar case it is possible to prove a
much stronger ill
posedness result (actually, a blow-up result).
\end{remark}

\begin{remark}\label{rem1}
   If in addition to \eqref{assN} we assume that the geodesic
   $\ga(s)$ is defined
   for all $s\in\R$, i.e., that the target manifold is complete (by the
   Hopf-Rinow Theorem), and that the dimension of the base space
   is $n=2$, then we can modify the proof of Theorem \ref{thm1}
   in such a way to use radial solutions exclusively.
   This is interesting in connection with the recent result of Struwe
   \cite{Str2}, who proved
   global existence of smooth radial solutions to the wave map system
   from $\R\times\R^{2}$ to the two dimensional sphere. Thus in this
   case the solution map is well defined for smooth radial data, but not
   uniformly continuous in $H^1$. (Actually, by a more complex construction,
   which involves a localization in
   Fourier space, it is possible to construct a radial counterexample
   also under assumption \eqref{assN} only). For a more
   precise statement we refer to the following proposition.
\end{remark}

\begin{prop}\label{prop1}
   Assume $n=2$ and the target space $N$ is complete.
   Then the conclusion of Theorem \ref{thm1} holds also
   if we restrict the solution map $\Phi$ to
   the subspace
   $H^{1}_{\mathrm{rad}}\times L^{2} _{\mathrm{rad}}$
   of radial functions in $H^{1}\times L^{2} $.
\end{prop}

\begin{remark}
It is important to notice that in the proof
of Theorem \ref{thm1} the fact that $\Gamma^{\ell}_{jk}$ are
Christoffel symbols of some Riemannian manifold is not essential.
In other words, the result holds for any system of the form
\eqref{wmsys},
provided the curves locally defined by the system of equations
$$\ga''_{\ell}+ \Gamma^\ell_{bc}(\ga)\ga'_{b}\ga'_{c}=0$$
satisfy an assumption like \eqref{gamma} near some point.
This means that the ill posedness in the sense of uniform continuity
is a general property of systems of the wave map
type.
\end{remark}

The authors are grateful to N.Tzvetkov for several discussions and
remarks concerning the well-posedness of semilinear problems.

\section{Proof of Theorem \ref{thm1}}

It is not restrictive to assume that $\ga$ is parameterized by arc length;
moreover, in the following we shall take $T=1$ for simplicity of
notations but the proof is unchanged in general. Assumption
~\eqref{gamma}
implies that for some component $\ga''_j$ of $\ga''$ one has
\begin{equation}\label{gamma1}
   |\ga''_j(s)|\geq c_1(N)\quad\hbox{for}\quad
      |s|\leq c_0(N)
\end{equation}
for suitable constants $c_0,c_1$ depending only on the manifold $N$.

Let $v,w$ be two $C^\infty$ real
valued solutions of the homogeneous wave equation
$$\square v=\square w=0$$
with data
$$v(0,x)=w(0,x)=v_0(x),\qquad \de_t v(0,x)=v_1(x),
    \qquad \de_t w(0,x)=w_1(x).$$
Notice that $v(0,x)\equiv w(0,x)$ and only the second datum is different.
Moreover, we shall always work with data of compact support, so
that $v(t,\cdot)$, $w(t,\cdot)$ will have support in a fixed
ball (say $B(0,10)$) for all $t\in[-1,1]$.
Then the functions $\uone=\ga\circ v$, $\utwo=\ga\circ w$
are solutions of the
wave map equation (see Remark \ref{remm1}), provided $v,w$ take their
values in the domain of $\ga(s)$; more precisely we shall assume
that
\begin{equation}\label{smallvw}
   |v|\leq c_{0}(N),\qquad |w|\leq c_{0}(N)
\end{equation}
and these conditions will be verified in the explicit construction
of $v$ and $w$. The corresponding Cauchy
data are given by
\begin{equation}\label{eqdata}
   \uone(0)=\utwo(0)=\ga(v_0),\quad \de_t\uone(0)=\ga'(v_0)v_1,\quad
      \de_t\utwo(0)=\ga'(v_0)w_1.
\end{equation}

Assume now that the solution map is defined
and uniformly continuous from some
neighbourhood
$U$ of 0 in $H^{n/2}(\R^n;N)\times H^{n/2-1}(\R^n;N)$ with values in
the space $C^1([0,1]; \dot H^{n/2-1})$
(the case of $C([0,1]; \dot H^{n/2})$ is
completely analogous).
If we apply this to the data
\eqref{eqdata}, we obtain that: \emph{for any $\ep>0$ there exists $\del>0$
such that}
\begin{equation}\label{unif}
    \|\de_t\uone(0)-\de_t\utwo(0)\|_{H^{n/2-1}}
       \leq\del\quad\Rightarrow\sup_{t\in[0,1]}
      \|\de_t\uone(t,\cdot)
           -\de_t\utwo(t,\cdot)\|_{\dot H^{n/2-1}}\leq\ep
\end{equation}
\emph{for all data $\uone(0)=\utwo(0)$ and $\de_t\uone(0)$, $\de_t\utwo(0)$ in $U$.}
We can express this condition in terms of the data for $v,w$. Indeed,
we have
$$\de_t\uone(0)-\de_t\utwo(0)=\ga'(v_{0})(v_{1}-w_{1}),$$
where $\ga'(v_{0})$ is a smooth function, equal to a constant outside
some compact set in $\R^{n}$. Applying Lemma \ref{commut}
in the Appendix for $s=n/2-1$, we have
\begin{equation*}
   \|\ga(f)g\|_{H^{n/2-1}}\leq c_{n}
      \| \ga(f)\|_{L^{\infty}\cap H^{n/2}}
          \cdot\|g\|_{H^{n/2-1}},
\end{equation*}
where we are using the notation
$$\|u\|_{X\cap Y}= \|u\|_{X}+ \|u\|_{Y}; $$
since $\ga(0)=0$ we can apply the standard Moser type estimate
\begin{equation}\label{MOS}
   \|\ga(f)\| _{H^{n/2}} \leq
     \rho_{0}(\|f\|_{L^{\infty}})\cdot\|f\|_{H^{n/2}}
\end{equation}
for a suitable continuous increasing function $\rho_{0}(s)$
(see e.g. \cite{Taylor97}, Vol.III, Chapter 13, Proposition 10.2),
we obtain an inequality like
\begin{equation}\label{GN}
   \|\ga(f)g\|_{H^{n/2-1}}\leq
   \rho_{1}(\|f\| _{L^{\infty}\cap H^{n/2}})\cdot\|g\|_{H^{n/2-1}}
\end{equation}
for some continuous increasing $\rho_{1}(s)$,
which is valid provided the range of the real valued function $f$
is contained in a compact subset of the domain of the smooth function
$\ga(s)$. Then we have
\begin{equation}\label{moser1}
   \|\de_t\uone(0)-\de_t\utwo(0)\|_{H^{n/2-1}}\leq
     \rho_{1}(\|v_{0}\| _{L^{\infty}\cap H^{n/2}})
        \|v_{1}-w_{1}\|_{H^{n/2-1}},
\end{equation}
hence property \eqref{unif} implies the following: for all $\ep>0$
there exists $\del>0$ such that
\begin{equation}\label{uniform}
   \|v_1-w_1\|_{H^{n/2-1}}\leq\del\quad\Rightarrow
   \sup_{t\in[0,1]}\|\de_t\uone(t,\cdot)
   -\de_t\utwo(t,\cdot)\|_{\dot H^{n/2-1}}\leq\ep
\end{equation}
for all data $(v_0,v_1)$ and $(w_0,w_1)$ belonging to a suitable
neighbourhood
$V$ of 0 in $H^{n/2}(\R^n)\times H^{n/2-1}(\R^n)$ and such
that $v_{0}=w_{0}$ and $|v_{0}|\leq c_{0}$, where
$c_{0}= c_{0}(N)$ is defined
in \eqref{gamma1}.

We now estimate from below the second term in \eqref{uniform}
$$\|\de_t\uone(t,\cdot)
   -\de_t\utwo(t,\cdot)\|_{\dot H^{n/2-1}} =
   \|\ga'(v)\de_{t}v
   -\ga'(w)\de_t w\|_{\dot H^{n/2-1}}. $$
We have
\begin{equation}\label{step1}
   \|\ga'(v)\de_{t}v
      -\ga'(v)\de_t w\|_{\dot H^{n/2-1}} \geq
         \|(\ga'(v)
      -\ga'(w))\de_t v\|_{\dot H^{n/2-1}} -
   \|\ga'(v)\de_{t}(v- w)\|_{\dot H^{n/2-1}}.
\end{equation}
We apply \eqref{GN} to the last term, obtainig
$$\|\ga'(v)\de_{t}(v-w)\|_{\dot H^{n/2-1}} \leq
  \rho_{1}(\|v\|_{L^{\infty}\cap H^{n/2}})
      \|\de_{t}(v-w)\|_{H^{n/2-1}};$$
by the energy identity for $\square(v-w)=0$,
$v_0=w_0$, we know that
$$\|\de_{t}(v-w)\|_{H^{n/2-1}} \leq
   \|v_{1}-w_{1}\|_{H^{n/2-1}} $$
and in conclusion
\begin{equation}\label{step2}
   \|\ga'(v)\de_{t}(v-w)\|_{\dot H^{n/2-1}} \leq
     \rho_{1}(\|v\| _{L^{\infty}\cap H^{n/2}})
        \|v_{1}-w_{1}\|_{H^{n/2-1}}
\end{equation}
To estimate from below the first term in the right side of
\eqref{step1} we use the Taylor developments
$$\ga'(b)-\ga'(a)=\ga''(a)(b-a)+F(a,b)(b-a)^{2},\qquad
  \ga''(a)=\ga''(0)+G(a)\cdot a$$
where $F(a,b), G(a)$ are smooth functions of their arguments whose
explicit expression is not relevant. Then
$$\ga'(v)-\ga'(w)=\ga''(0) \cdot(v-w) +R(v,w)\cdot(v-w)$$
where we have written for short
$$R(u,v)= G(v)\cdot v+F(v,w)(v-w) .$$
Recalling that $|\ga''(0)|\geq c_{1}$
(see \eqref{gamma1}), we have
$$\|\ga''(v)(v-w)\de_{t}v\| _{\dot H^{n/2-1}}\geq
   c_{1} \|(v-w)\de_{t}v\| _{\dot H^{n/2-1}}-
   \| R(v,w)(v-w)\de_{t}v\| _{\dot H^{n/2-1}}.
   $$
Now we can apply
\eqref{multest} of Lemma \ref{commut} in the Appendix to obtain
$$\|R(v,w)(v-w)\de_{t}v\| _{\dot H^{n/2-1}}\leq
   \|R(v,w)\|_{L^{\infty}\cap H^{n/2}}\|(v-w)\de_{t}v\| _{H^{n/2-1}}
$$
while using \eqref{MOS} it is standard to obtain
$$\|R(v,w)\|_{L^{\infty}\cap H^{n/2}} \leq
   \rho_{2}(\| v,w\| _{L^{\infty}\cap H^{n/2}})\cdot
   \| v,w\| _{L^{\infty}\cap H^{n/2}}$$
($\|v,w\|=\|v\|+\|w\|$)
for some continuous increasing function $\rho_{2}(s)$
whose precise form is not
relevant. In conclusion, recalling also \eqref{step1} and
\eqref{step2}, we have proved the inequality
\begin{multline}\label{goodstep}
  \|\de_t\uone(t,\cdot)
   -\de_t\utwo(t,\cdot)\|_{\dot H^{n/2-1}}\geq
   c_{1} \|(v-w)\de_{t}v\| _{\dot H^{n/2-1}}-\\
      -\rho_{2}(\| v,w\| _{L^{\infty}\cap H^{n/2}})\cdot
   \| v,w\| _{L^{\infty}\cap H^{n/2}}
   \|(v-w)\de_{t}v\| _{H^{n/2-1}}-\\
     -\rho_{1}(\|v\| _{L^{\infty}\cap H^{n/2}})
        \|v_{1}-w_{1}\|_{H^{n/2-1}}.
\end{multline}
To proceed, we must construct explicitly the functions $v$ and $w$.
This is done with the help of a few lemmas.

\begin{lemma}\label{lem1}
   Let $n\geq2$. There exists a sequence of real valued
   functions
   $\phi_j\in C^\infty_0(\R^n)$ supported
   in the ball $\{|x|\leq 2\}$, with
   \begin{equation}\label{phito0}
     \phi_j\to0\quad\hbox{in}\quad H^{n/2-1}(\R^n)\quad
     \hbox{as}\quad j\to\infty
   \end{equation}
   such that,
   denoting by $z_j(t,x):\R\times\R^n\to\R$ the solution of
   the linear problem
   \begin{equation}\label{homeq}
     \square z=0,\qquad z(0,x)=0,\qquad \de_t z(0,x)=\phi_j(x)
   \end{equation}
   one has
   \begin{equation}\label{zeq1}
     z_j(1,0)=1\qquad\hbox{for any}\quad j.
   \end{equation}
   The functions $\phi_j$ and hence $z_j(t,x)$ can be chosen as radial
   functions in $x$, i.e., depending only on $|x|$.
\end{lemma}

\begin{proof}
   We begin by the case $n=2$. For $0<p<q<1$, we define $\psi_{p,q}(y)$
   on $\R^2$ as follows:
   \begin{equation}\label{psipq}
     \psi_{p,q}(y)=-\frac{I_{\{p\leq|y|\leq q\}}(y)}
            {\sqrt{1-|y|^2}\log(1-|y|^2)}
   \end{equation}
   where $I_A(y)$ denotes the characteristic function of the set $A$.
   An elementary computation gives
   \begin{equation}\label{L2norm}
     2\|\psi_{p,q}(x)\|^2_{L^2(\R^2)}=
        \frac1{\log{(1-q^2)}}-\frac1{\log{(1-p^2)}}.
   \end{equation}
   Notice that taking any $0<p_1<p<q<q_1<1$,
   and an arbitrary smooth radial cutoff function $\chi_{p,q}$ with
   $$I_{\{p\leq|y|\leq q\}}(y)\leq\chi_{p,q}(y)
       \leq I_{\{p'\leq|y|\leq q'\}}(y),$$
   we can modify definition ~\eqref{psipq}
   as follows:
   \begin{equation}\label{psitilde}
     \widetilde\psi_{p,q}(y)=-\frac{\chi_{\{p\leq|y|\leq q\}}(y)}
            {\sqrt{1-|y|^2}\log(1-|y|^2)}
   \end{equation}
   in order to obtain a smooth initial datum with similar norm:
   \begin{equation}\label{L2normtil}
     \frac1{\log{(1-q^2)}}-\frac1{\log{(1-p^2)}}
       \leq 2\|\widetilde\psi_{p,q}(x)\|^2_{L^2(\R^2)}
       \leq\frac1{\log{(1-q_1^2)}}-\frac1{\log{(1-p_1^2)}}.
   \end{equation}

   On the other hand, the solution $z_{p,q}(t,x)$ of the problem
   $$\square z=0,\qquad z(0,x)=0,\qquad \de_t z(0,x)=\psi_{p,q}$$
   is explicitly given by
   $$z_{p,q}(t,x)=\frac t{2\pi}\int_{|x-y|\leq t}
       \frac{\psi_{p,q}(y)}{\sqrt{t^2-|x-y|^2}}dy$$
   and in particular at $(t,x)=(1,0)$ one has
   $$z_{p,q}(1,0)=-\frac 1{2\pi}\int_{p\leq|y|\leq q}
      \frac{1}{(1-|y|^2)\log(1-|y|^2)}dy=
      \frac 1{4\pi}
      \log\left|\frac{\log(1-q^2)}{\log(1-p^2)}\right|.$$
   By the positivity of the kernel we have immediately, for the solution
   $\widetilde z_{p,q}$ obtained by replacing $\psi_{p,q}$
   with $\widetilde\psi_{p,q}$,
   $$\frac 1{4\pi}
      \log\left|\frac{\log(1-q^2)}{\log(1-p^2)}\right|
      \leq \widetilde z_{p,q}(1,0) \leq
      \frac 1{4\pi}
      \log\left|\frac{\log(1-q_1^2)}{\log(1-p_1^2)}\right|. $$
   If we now choose for $\del\in]0,1[$
   $$1-p_1^2=\del, \qquad 1-p^2=\del^2, \qquad 1-q^2=\del^3,
       \qquad 1-q_1^2=\del^4$$
   and write $\psi_\del=\widetilde\psi_{p,q}$, we obtain
   $$\frac1{\sqrt{12}}|\log\del|^{-1/2}\leq
       \|\widetilde\psi_\del\|_{L^2}\leq
       \sqrt{\frac38}|\log\del|^{-1/2}\to0\qquad
      \hbox{as }\del\to0$$
   while $z_\del=\widetilde z_{p,q}$ satisfies
   $$\frac1{4\pi}\log\frac32\leq z_\del(1,0)\leq\frac1{4\pi}\log4.$$
   Defining $\phi_j=z_{\del_j}(1,0)^{-1}\psi_{\del_j}$ for any
   $\del_j\downarrow0$
   we obtain the thesis.

   The general case for $n$ even, $n\geq2$
   follows easily by modifying the above
   example, using the fact that the solution of~\eqref{homeq} can be
   represented as
   $$z(t,x)=\sum_{0\leq |\al|\leq (n-2)/2}a_\al t^{|\al|+1}
       \int_{|y|\leq1} y^\al D^\al\phi(x+ty)(1-|y|^2)^{-1/2}dy,$$
   which, for a radial function $\phi$, gives
   $$z(1,0)=\sum_{\nu=0}^{(n-2)/2}c_{\nu}
       \int_0^1\de_r^\nu\phi(r)(1-r^2)^{-1/2}r^{\nu+n-1}dr.$$
   Here of course we shall
   choose a datum $\phi$ such that its radial derivative of order
   $n/2-1$ is
   of the form $\widetilde\psi_{p,q}$ seen above.

   Let us now consider the case of odd $n$, starting from $n=3$.
   In this case it is sufficient to use the well known fact (see
   e.g., Theorem 11.1 in Volume I of \cite{LiMa}) that,
   for any bounded $\Ome\subset\R^n$ with $C^\infty$
   boundary,
   $$C^\infty_0(\Ome)\quad\hbox{is dense in}\quad H^{1/2}(\Ome)$$
   and also in $H^{s}$ for $s\leq1/2$. Since we shall need a special
   version of this result for radial functions, we shall give here
   a self-contained proof adapted to our situation.

   Indeed, consider the space
   \begin{equation}
      Z=\{\phi\in C_{0}^{\infty}(\R)\;:\;\phi(x)=\phi(-x),\
          \phi\equiv0\hbox{\ near\ }1\hbox{\ and\ }-1\}
   \end{equation}
   (where ``near $\pm1$'' means
   ``on some neighbourhood of these two
   points, depending on $\phi$''). It is easy to see that $Z$ is a dense
   subset of the space of
   even $H^{1/2}(\R)$ functions
   \begin{equation}
      H^{1/2}_{\mathrm{even}}(\R)=
      \{u\in H^{1/2}(\R)\;:\;u(x)=u(-x)\}
   \end{equation}
   by the following argument: in the Hilbert space
   $H^{1/2}_{\mathrm{even}}(\R)$ we can certainly choose a $u_{0}$
   orthogonal to $Z$, and we must only prove that $u_{0}=0$. The
   tempered distribution $T$ whose Fourier transform is given by
   $$\widehat T=\langle\xi\rangle\widehat u_{0}$$
   belongs to $H^{-1/2}(\R)$ and by the identity
   $$T(\phi)=(\langle\xi\rangle\widehat u_{0},
                            \widehat{\overline{\phi}})_{L^{2}}=
        (u_{0},\overline\phi)_{H^{1/2}}=0$$
   for any test function in $Z$, we see that the support of $T$ is
   contained in the set $\{\pm1\}$, i.e., $T$ is a linear combination of
   a finite number of derivatives of $\del_{1}$, $\del_{-1}$. Hence
   $\widehat T(\xi)$ is a function of the form
   $$\widehat T(\xi)=\sum_{\ell=0}^{N}
      (c_{\ell}e^{i\xi}+ d_{\ell}e^{-i\xi})\xi^{\ell}$$
   for a suitable $N\geq0$ and complex numbers $c_{\ell},d_{\ell}$,
   and at the same time
   $\langle\xi\rangle^{-1/2}\widehat T(\xi)$ must belong to $L^{2}$.
   It is trivial to see that the only such function $\widehat T$
   is 0, and this implies $u_{0}\equiv0$ too.

   Thus we have proved that $C_{0}^{\infty}(]1,1[)$ is dense in
   $H^{1/2}_{\mathrm{even}}(]-1,1[)$ since this last space coincides with
   the space of restrictions of functions in $H^{1/2}_{\mathrm{even}}(\R)$
   to $]-1,1[$, with the restriction norm (the norm of $u$ is the
   infimum
   of the norms of its possible extensions).

   It would not be difficult to prove the same result
   for higher dimensions, but actually here we only need to construct
   a sequence of radial functions $\psi_{j}\in C^{\infty}_{0}(B_{1})$
   which converges to 1 in $H^{1/2}(B_{1})$, where $B_{1}$ is the
   unit ball $ B_{1}=\{x\in\R^{n}:|x|<1\}$.    To this end, it
   is sufficient to remark that the operator
   $$A:H^{s}_{\mathrm{even}}(]-1,1[)\to H^{s}(B_{1})$$
   defined as
   $$A(f)(x)=f(|x|)$$
   is bounded for all $0\leq s\leq 1$: this is proved directly for
   $s=0,1$ and follows e.g. by interpolation for the intermediate
   values of
   $s$. Hence taken any sequence $f_{j}(x)$ in
   $C^{\infty}_{0}(]-1,1[)$, with $f(x)=f(-x)$, converging to 1
   in the $H^{1/2}(]-1,1[)$ norm, we need only define
   $$\psi_{j}(x)=f_{j}(|x|)$$ to obtain the desired result.

   Now, setting $\phi_j=1-\psi_j$ we obtain
   a sequence of radial smooth functions on $B_1$, converging to 0 in
   the norm of
   $H^{1/2}(B_1)$, and identically equal to 1 on some neighbourhood
   of $\de B_1$
   (depending on $j$). By Kirchhoff's formula then we obtain
   $$z(1,0)=\frac1{4\pi}\int_{\de B_1} \phi_j(y)dS=1$$
   as needed.

   In the general case $n\geq3$ odd, we proceed in a similar way
   using the
   general representation of the solution; notice that for radial
   $\phi$ the following formula holds
   $$z(1,0)=\sum_{\nu=0}^{(n-3)/2}b_\nu\de_r^\nu\phi(1)$$
   for suitable constants $b_{\nu}$.
\end{proof}

\goodbreak

In the construction of the Lemma we have no control on the $L^{\infty}$
norm of the functions $z_{j}$; if we give up the requirement that the
$z_{j}$ be radial, however, it is easy to obtain the following result:

\begin{cor}\label{corol1}
   Let $n\geq2$. There exists a sequence of real valued
   functions
   $\phi_j\in C^\infty_0(\R^n)$ supported in the ball $\{|x|\leq5\}$
   with
   \begin{equation}\label{phito0bis}
     \phi_j\to0\quad\hbox{in}\quad H^{n/2-1}(\R^n)\quad
     \hbox{as}\quad j\to\infty
   \end{equation}
   such that,
   denoting by $z_j(t,x):\R\times\R^n\to\R$ the solution of
   the linear problem
   \begin{equation}\label{homeqbis}
     \square z=0,\qquad z(0,x)=0,\qquad \de_t z(0,x)=\phi_j(x)
   \end{equation}
   one has
   \begin{equation}\label{zeq1bis}
     z_j(t_{j},0)=1\qquad\hbox{for some sequence }\quad t_{j}\in]0,1]
   \end{equation}
   and
   \begin{equation}\label{zeq2bis}
      |z_{j}(t,x)|\leq 1\qquad \hbox{for all }
         (t,x,j)\in[0,1]\times\R^{n}\times\N.
   \end{equation}
\end{cor}

\begin{proof}
  The functions $z_{j}$ constructed in the Lemma are smooth and
  compactly supported, let $(t_{j},x_{j})$ be a point where $|z_{j}|$
  attains its maximum value $m_{j}$ on the strip
  $[0,1]\times\R^{n}$, and define
  $$\widetilde{z_{j}}(t,x)=m_{j}^{-1}z_{j}(t,x-x_{j})$$
  (and possibly multiply by the sign of
  $z_{j}(t_{j},x_{j})$).
  Notice that $t_{j}>0$ since $z(0,x)\equiv0$.
  This concludes the proof.
\end{proof}

Before passing to the main body of the proof, a last
elementary rescaling lemma is necessary.

\begin{lemma}\label{lem2}
   Let $\chi(x)$ with $\|\chi\|_{\dot H^{n/2-1}}\neq0$
   be a smooth compactly supported (radial)
   function, vanishing for $|x|\geq2$, and with the property
   \begin{equation}\label{chiz}
      \int_{\R^{n}}\chi(x)dx=0.
   \end{equation}
   Let $R\geq1$, $M\geq0$ be positive numbers,
   $0\leq T\leq1$, and denote by $\vrm(t,x)$
   the (radial) solution
   of the homogeneous wave equation
   \begin{equation}\label{reverse}
      \square v=0,\qquad v(T,x)=0,\qquad
      \de_{t}v(T,x)=\chi_{R,M}(x)\equiv M\cdot\chi(Rx)
   \end{equation}
   with data at $t=T>0$. Denote by $v_{0},v_{1}$ the traces
   \begin{equation}\label{datav}
      v_0=\vrm(0,x),\qquad v_1=\de_t\vrm(0,x)
   \end{equation}
   so that \eqref{reverse} is equivalent to a Cauchy problem for
   the homogeneous wave equation with data $v_{0},v_{1}$ at $t=0$.
   Then the following estimates hold, for a constant $c_{n}$
   depending only on the space dimension $n$ and on the
   function $\chi(x)$:
   \begin{equation}\label{init}
      \|v_{0}\|_{H^{n/2}}+\|v_{1}\|_{H^{n/2-1}}\leq
         c_{n}\frac M {R},
   \end{equation}
   and, for all $(t,x)\in[0,1]\times\R^{n}$,
   \begin{equation}\label{estLinf}
      |\vrm(t,x)|\leq c_{n}\frac M {R}.
   \end{equation}
   Finally, for all $0\leq s\leq n/2$ and all $t\in\R$
   \begin{equation}\label{alls}
      \|v(t,\cdot)\|_{\dot H^{s}}+
      \|\de_{t}v(t,\cdot)\|_{\dot H^{s-1}}\leq c_{n}
      \frac M R \cdot R^{s-n/2}.
   \end{equation}
\end{lemma}

\begin{proof}
   Rescale $v(t,x)$ as
   $$v(t,x)=w(Rt,Rx)$$
   so that
   $$\square w=0,\qquad w(RT,x)=0,\qquad
     \de_{t}w(RT,x)=\frac M R\chi(x).$$
   By the energy estimates we have for all real $s$ and all
   $t\in\R$
   \begin{equation}\label{nrg}
      \|w(t,\cdot)\|_{\dot H^{s}}+
      \|\de_{t}w(t,\cdot)\|_{\dot H^{s-1}}\leq
             2\frac M R\|\chi\|_{\dot H^{s-1}}
   \end{equation}
   which scaling back to $v$ gives
   $$\|v(t,\cdot)\|_{\dot H^{s}}+
     \|\de_{t}v(t,\cdot)\|_{\dot H^{s-1}}\leq
     2\frac M R\cdot R^{s-n/2}\|\chi\|_{\dot H^{s-1}}.$$
   Notice that \eqref{nrg} gives a finite bound
   also for $s=0$; indeed, by
   assumption \eqref{chiz} we have $\widehat\chi(0)=0$ and hence
   $\widehat\chi/|\xi|\in L^{2}$, i.e., $\chi\in\dot H^{-1}$.
   Thus for all $0\leq s\leq n/2$ we obtain ($R\geq1$)
   $$\|v(t,\cdot)\|_{\dot H^{s}}+
     \|\de_{t}v(t,\cdot)\|_{\dot H^{s-1}}\leq
     2\frac M R R^{s-n/2}\cdot
     (\|\chi\|_{\dot  H^{n/2-1}}+\|\chi\|_{\dot H^{-1}}). $$
   This proves \eqref{alls}; inequality \eqref{init}
   is just the special case $s=n/2$ computed at $t=0$.

   To prove \eqref{estLinf} we use \eqref{nrg} again
   for $s=n/2+1$, which gives
   $$\sup_{t\in\R}\|w(t,\cdot)\|_{\dot H^{n/2+1}}\leq
     c_{n}\frac M R\|\chi\|_{\dot H^{n/2}},$$
   while for $s=0$ it gives
   $$\sup_{t\in\R}\|w(t,\cdot)\|_{L^{2}}\leq
     c_{n}\frac M R\|\chi\|_{\dot H^{-1}},$$
   and this is bounded by \eqref{chiz} as already remarked.
   Thus, by Sobolev embedding, we have
   $$\|w\|_{L^{\infty}(\R\times\R^{n})}\leq
          c_{n}\sup_{t\in\R}\|w(t,\cdot)\|_{H^{n/2+1}}\leq
          c'_{n}\frac M R(\|\chi\|_{\dot H^{n/2}}
          +\|\chi\|_{\dot H^{-1}}).$$
   Since $\|v\|_{L^{\infty}}=\|w\|_{L^{\infty}}$, this concludes
   the proof; the constant $c_{n}$ depends only on $n$ and the quantity
   $\|\chi\|_{\dot  H^{n/2}}+\|\chi\|_{\dot H^{-1}}$.
\end{proof}

We revert now to the main proof. The next step is
the explicit construction
of sequences of functions $v,w$ appearing in
\eqref{goodstep}. As $v$ we shall choose the function $\vrm$ constructed
in the preceding lemma, with a suitable choice of the parameters.
Notice that by \eqref{init}
we can assume that the initial data $v_0,v_1$ belong to the
neighbourhood $V$ of 0 in $H^{n/2}\times H^{n/2-1}$ on which property \eqref{uniform}
holds, as soon as $M/R$ is small enough; e.g., if $V$ contains a ball
of radius $r_0(V)$ centered in 0, we may assume that
\begin{equation}\label{datal2}
   4c_{n}\|\chi\|_{H^{n/2-1}}\frac M R< r_0.
\end{equation}
Notice also that, in order to define the composition $\ga\circ v$,
we must ensure that $|v|<s_{0}$ (at least on the strip
$[0,1]\times\R^{n}$) since the geodesic curve is only defined
on the interval $]-s_{0},s_{0}[$, or even better, that $|v|<c_{0}$
given by \eqref{gamma1}.
Using \eqref{estLinf}, we see that it is sufficient to further decrease
$M/R$, e.g., to impose the condition
\begin{equation}\label{datalinf}
   c_{n}\|\chi\|_{H^{n/2}}\frac M R< c_0/2.
\end{equation}
In connection with Remark \ref{rem1}, we observe that condition
\eqref{datalinf} is not necessary when we assume that $\ga(s)$ is defined
for all $s\in\R$.

Then we define $v_{j}=\vrm$ with the following choices.
The parameter $T$ will be chosen as
\begin{equation}\label{choiceT}
   T=t_{j}\qquad\Rightarrow \qquad \vrm(t_{j},x)=M\chi(Rx),
\end{equation}
where $t_{j}$ are given by Corollary \ref{corol1};
the parameter $R=R_{j}$ will be chosen such that
\begin{equation}\label{choiceR}
   z_{j}(t_{j},x)\geq\frac12\hbox{\ \ on the ball\ \ }
      \{|x|\leq2{R_{j}}^{-1}\};
\end{equation}
this is possible in view of \eqref{zeq1bis}
and of the continuity of $z_{j}$; it is
not restrictive to assume that $R_{j}\uparrow+\infty$.
On the parameter $M=M_{j}$, besides \eqref{datal2}, \eqref{datalinf}
further smallness conditions will be imposed in the following.

We now define $w_{j}$; let $\mu>0$ be a small parameter, and set
($w_0\equiv v_0$ and)
\begin{equation}\label{choice}
   w_j=v_j+\mu z_j
\end{equation}
where $v_{j}$ was defined above and $z_{j}$ is given by Corollary
\ref{corol1}.
Thus the data for $w_{j}$ are
$$w_{0}\equiv v_{0},\quad
  w_{1}=v_{1}+\mu\phi_{j}$$
with $v_{0}, v_{1}$ the traces of $v_{j}$ at $t=0$, studied in
Lemma \ref{lem2}. Again, in order to define
the composition $\ga\circ w$, we must ensure that $|w|<c_{0}$, at
least for $0\leq t\leq1$. Using
\eqref{zeq2bis} and recalling
\eqref{estLinf}, \eqref{datalinf}, we see that it is
sufficient to impose the condition
\begin{equation}\label{mu}
   0<\mu<c_{0}/2.
\end{equation}
Notice that the data $w_{0},w_{1}$ belong to the given neighbourhood
$V$ as soon as $j$ is large enough, since $\phi_{j}\to0$ in
$H^{n/2-1}$.

Consider inequality \eqref{goodstep}; our aim is to estimate
its right side from below. The first term at $t=t_{j}$ gives
\begin{equation}\label{1stt}
   \|(v_{j}-w_{j})\de_{t}v_{j}\|_{\dot H^{n/2-1}}=
   \mu \|z_{j}(t_{j},\cdot)\chi_{R_{j},M_{j}}\| _{\dot H^{n/2-1}};
\end{equation}
we can apply \eqref{below2} of the Appendix, with $s=n/2-1$;
since $z_{j}\geq 1/2$ on the support of $\chi_{R,M}$, we have
$$\|z_{j}(t_{j},\cdot)\chi_{R _{j},M _{j}}\| _{\dot H^{n/2-1}}
  \geq\frac c2\|\chi_{R _{j},M _{j}}\| _{\dot H^{n/2-1}}
      -c'\|z_{j}\|_{H^{n/2}}\|\chi_{R _{j},M _{j}}\| _{H^{n/2-1}} .$$
Now we have for $R$ large enough
$$\|\chi_{R,M}\| _{\dot H^{n/2-1}} =\frac MR\kappa,\qquad
  \|\chi_{R,M}\| _{H^{n/2-1}} \leq 2 \frac MR\kappa,$$
where by assumption
$$\kappa= \|\chi\| _{\dot H^{n/2-1}} \neq0.$$
Moreover, by the energy identity we have for all $t$
\begin{equation}\label{energz}
   \|z_{j}(t,\cdot)\|_{\dot H^{n/2}}\leq c\|\phi_{j}\|_{\dot H^{n/2-1}}
   \to 0\hbox{\ \ as\ \ }j\to \infty
\end{equation}
and also for all $|t|\leq1$
\begin{equation}\label{l2z}
   \|z_{j}(t,\cdot)\|_{L^{2}}=
   \left\|\frac{\sin(t|\xi|)}{|\xi|}
         \widehat\phi_{j}\right\| _{L^{2}}\leq
   \| \phi_{j}\| _{L^{2}}\to0\qquad\Rightarrow\qquad
   \|z_{j}(t,\cdot)\|_{H^{n/2}}\to0
\end{equation}
Hence we have proved that
\begin{equation}\label{2stt}
   \|(v_{j}-w_{j})\de_{t}v_{j}\|_{\dot H^{n/2-1}}\geq
    c (n,\kappa) \frac{M_{j}}{R_{j}}.
\end{equation}

In view of the second term in \eqref{goodstep} we need also a bound
from above for the quantity $(v-w)\de_{t}v$; by \eqref{multest} from the
Appendix, with $s=n/2-1$, we have
$$\|z_{j}(t_{j},\cdot)\chi_{R_{j},M_{j}}\| _{H^{n/2-1}}
   \leq C\|z_{j}\|_{L^{\infty}\cap H^{n/2}}
   \|\chi_{R_{j},M_{j}}\|_{H^{n/2-1}}
   \leq C\frac{M_{j}}{R_{j}} \|z_{j}\|_{L^{\infty}\cap H^{n/2}} ,$$
for $j$ large enough, and recalling that $z_{j}\to0$ in
$H^{n/2}$ uniformly in $|t|\leq1$ as remarked above, and
$|z_{j}|\leq1$ by construction, we finally obtain
\begin{equation}\label{3stt}
   \|(v_{j}-w_{j})\de_{t}v_{j}\|_{H^{n/2-1}}\geq
    c' (n,\kappa) \frac{M_{j}}{R_{j}}
\end{equation}
provided $j$ is large enough.

We notice that, by \eqref{estLinf}, \eqref{alls},
$$\|v_{j}\|_{L^{\infty}\cap H^{n/2}}\leq c\frac {M_{j}}{R_{j}} $$
while, recalling that $|z_{j}|\leq1$ and that $\|z_{j}\|_{H^{n/2}}\leq1$
for $j$ large enough, we have
$$\|w_{j}\| _{L^{\infty}\cap H^{n/2}} =
  \|v_{j}+\mu z_{j}\| _{L^{\infty}\cap H^{n/2}} \leq
  c\mu + c\frac {M_{j}}{R_{j}} .$$
Together with \eqref{3stt} this gives us the following estimate for the
second term in \eqref{goodstep}:
\begin{multline}\label{4stt}
   \rho_{2}(\| v_{j},w _{j}\| _{L^{\infty}\cap H^{n/2}})
   \| v _{j},w _{j}\| _{L^{\infty}\cap H^{n/2}}
   \|(v_{j}-w_{j})\de_{t}v_{j}\|_{H^{n/2-1}}\leq\\ \leq
   \rho_{3}(\mu+M_{j}/R_{j})\cdot
      \left(\mu+ \frac{M_{j}}{R_{j}}\right) \frac{M_{j}}{R_{j}}.
\end{multline}
We can impose now the last smallness condition on $\mu$ and
$M_{j}$ (recall that $M_{j}/R_{j}$ is bounded):
\begin{equation}\label{lastcond}
   \rho_{3} (\mu+M_{j}/R_{j})\cdot
      \left(\mu+ \frac{M_{j}}{R_{j}}\right)
      \leq \frac12 c(n,\kappa)
\end{equation}
where $c(n,\kappa) $ is the constant appearing in
\eqref{2stt}. Thus we get
\begin{equation}\label{5stt}
   \rho_{2}(\| v_{j},w _{j}\| _{L^{\infty}\cap H^{n/2}})
   \| v _{j},w _{j}\| _{L^{\infty}\cap H^{n/2}}
  \|(v_{j}-w_{j})\de_{t}v_{j}\|_{H^{n/2-1}}\leq
        \frac12 c(n,\kappa)
      \frac{M_{j}}{R_{j}}.
\end{equation}

The last term in \eqref{goodstep} is quite easy to estimate:
we have for $j\to\infty$
\begin{equation}\label{6stt}
   \rho_{1}(\|v\| _{L^{\infty}\cap H^{n/2}})
        \|v_{1}-w_{1}\|_{H^{n/2-1}}\leq
     \rho_{4}(\mu+M_{j}/R_{j})\cdot\mu\|\phi_{j}\| _{H^{n/2-1}}\to0.
\end{equation}

We can finally choose $M_{j}$ and $\mu$;
we define $M_{j}=\la\cdot R_{j}$, and $\la$, $\mu$ are two
positive constants so small that conditions
\eqref{datal2}, \eqref{datalinf}, \eqref{lastcond} are
satisfied.

Summing up, by \eqref{2stt}, \eqref{5stt}, \eqref{6stt}, we
obtain
\begin{equation}\label{lllast}
   \|\de_t\uone(t_{j},\cdot)
   -\de_t\utwo(t_{j},\cdot)\|_{\dot H^{n/2-1}}\geq
   \frac14 c(n,\kappa)
      \frac{M_{j}}{R_{j}}= \frac14 c(n,\kappa)\cdot\la.
\end{equation}
provided $j$ is large enough.

We can now conclude the proof. Recalling
\eqref{datav}, we can choose as data for $v$ the sequences
$$v_0^{(j)}=v_{R_j,M_j,t_{j}}(0,x),\qquad
   v_1^{(j)}=\de_t v_{R_j,M_j,t_{j}}(0,x)$$
while the data for $w$ are chosen as
$$w_0^{(j)}=v_0^{(j)},\qquad
   w_1^{(j)}=v_1^{(j)}+\mu\phi_j    \equiv v_1^{(j)}+c_0\phi_j/2.$$
By \eqref{datal2} the data for $v$ belong to $V$; as
a consequence, the
data for $w$ belong to $V$ provided $j$ is large enough, since
$\phi_j\to0$ in $H^{n/2-1}$.
Thus we are in position to apply the uniform continuity property
~\eqref{uniform}; since $w_1-v_1=\mu\phi_j/2$ we have that
for all $\ep>0$ there exists $\del>0$ such that
$$\|\phi_j\|_{H^{n/2-1}}<\del\qquad\Rightarrow\qquad
    \sup_{t\in[0,1]}\|\de_t\uone(t,\cdot)-
    \de_t\utwo(t,\cdot)\|_{\dot H^{n/2-1}}\leq\ep;$$
hence in particular at $t=t_{j}$ we must have
$$\|\phi_j\|_{H^{n/2-1}}<\del\qquad\Rightarrow\qquad
    \|\de_t\uone(t_{j},\cdot)
        -\de_t\utwo(t_{j},\cdot)\|_{\dot H^{n/2-1}}\leq\ep$$
and this is in clear contradiction with ~\eqref{lllast}.

\section{Proof of Proposition \ref{prop1}}

The proof follows exactly the same lines as for Theorem \ref{thm1},
and actually it is simpler from a technical point of view.
Indeed, when $n=2$ we must violate the following uniform continuity
condition:
\emph{for any $\ep>0$ there exists $\del>0$
such that}
$$\|\de_t\uone(0)-\de_t\utwo(0)\|_{L^2}\leq\del\quad\Rightarrow
\sup_{t\in[0,1]}\|\de_t\uone(t,\cdot)-\de_t\utwo(t,\cdot)\|_{L^2}\leq\ep$$
\emph{for all data $\uone(0)=\utwo(0)$ and $\de_t\uone(0)$, $\de_t\utwo(0)$
in $U$.}
We choose as above two $C^\infty$ real
valued solutions $v,w$ of the homogeneous wave equation
$$\square v=\square w=0$$
with data
$$v(0,x)=w(0,x)=v_0(x),\qquad \de_t v(0,x)=v_1(x),
    \qquad \de_t w(0,x)=w_1(x),$$
and we set $\uone=\ga\circ v$ , $\utwo=\ga\circ w$. Since
$|\ga'(s)|=1$ by the choice of the arc lenght parameter,
the uniform continuity for our choice of data
becomes simply
\begin{equation}\label{uniform2}
   \|v_1-w_1\|_{L^2}\leq\del\quad\Rightarrow
   \sup_{t\in[0,1]}\|\de_t\uone(t,\cdot)-
        \de_t\utwo(t,\cdot)\|_{L^2}\leq\ep
\end{equation}
In order to violate this property,
we estimate from below the second term in \eqref{uniform2}.
We can write
$$\de_t\uone-\de_t\utwo=\ga'(v)v_t-\ga'(w)w_t=
   (\ga'(v)-\ga'(w))w_t+\ga'(w)\cdot\de_t(v-w)$$
whence
$$|\de_t\uone-\de_t\utwo|\geq|\ga'(v)-\ga'(w)|\cdot|w_t|-|\de_t(v-w)|$$
using the identity $|\ga'|\equiv1$. This implies easily
$$|\de_t\uone-\de_t\utwo|^2+|\de_t(v-w)|^2
       \geq\frac12|\ga'(v)-\ga'(w)|^2\cdot|v_t|^2,$$
which can be written
\begin{equation}\label{below}
   |\de_t\uone-\de_t\utwo|^2+|\de_t(v-w)|^2
       \geq\frac12\left|\int_v^w\ga''(\sig)d\sig\right|^2\cdot|v_t|^2.
\end{equation}

As $v$ we shall choose the radial function
$\vrm$ constructed
in Lemma \ref{lem2}, with a suitable choice of the parameters.
By \eqref{init} for $n=2$
we can assume that the initial data $v_0,v_1$ belong to the
neighbourhood $V$ of 0 in $H^1\times L^2$ on which property \eqref{uniform}
holds, as soon as $M/R$ is small enough; if $V$ contains a ball
of radius $r_0(V)$ centered in 0, we may assume that
\begin{equation}\label{datal22}
   4\|\chi\|_{L^{2}}\frac M R< r_0.
\end{equation}
Notice that, thanks to the assumption that the geodesic curve
is globally defined, it is not necessary to impose any restriction
to the $L^{\infty}$ norm of $v$.

We now choose the data for $w$; let $\mu>0$ be a small parameter, and set ($w_0\equiv v_0$ and)
\begin{equation}\label{choice2}
   w_1=v_1+\mu\phi_j
\end{equation}
where $\phi_j$ are the smooth radial functions constructed in Lemma
\ref{lem1}.
Then we have $w=v+\mu z_j$, with $v\equiv\vrm$; again, no condition
on the $L^{\infty}$ norm of $w$ is necessary since $\ga(s)$
is defined for all $s$.

Recall now \eqref{below} which gives
\begin{equation}\label{below22}
   |\de_t\uone-\de_t\utwo|^2+|\mu\de_t z_j|^2
       \geq\frac12\left|\int_v^{v+\mu z_j}
       \ga''(\sig)d\sig\right|^2\cdot|v_t|^2.
\end{equation}
Notice that this is a pointwise inequality, valid at any $(t,x)$.
We can fix now $t=1$, choose $T=1$ in the definition of
$\vrm$ while leaving $R,M$ free for the moment (apart from
\eqref{datal22}), and we get
\begin{equation}\label{below3}
   |\de_t(\uone-\utwo)(1,x)|^2+|\mu\de_t z_j(1,x)|^2
       \geq\frac12\left|\int_0^{\mu z_j(1,x)}
       \ga''(\sig)d\sig\right|^2\cdot|\chi_{R,M}|^2.
\end{equation}
Integrating on $\R^2$, and using the energy inequality
$\|\de_t z_j(t,\cdot)\|_{L^2}\leq\|\phi_j\|_{L^2}$ we get
($\|\cdot\|=\|\cdot\|_{L^2}$)
\begin{equation}\label{below4}
   \|\de_t(\uone-\utwo)(1,\cdot)\|^2+\|\mu\phi_j\|^2
       \geq\frac12\left\|
       \int_{0}^{\mu z_j(1,\cdot)}
                 \ga''(\sig)d\sig\cdot \chi_{R,M}\right\|^2.
\end{equation}
Since $z_j$ is smooth and satisfies \eqref{zeq1}, we also have
\begin{equation}\label{zeqx}
   2> z_j(t_{j},x)>1/2\hbox{\ \ for\ \ }|x|\leq \frac1{R_j}
\end{equation}
for some $R_j$ large enough;
this is our choice for the parameter $R=R_{j}$ in the definition
of $\vrm$. Moreover, we shall choose $M=M_j$ proportional
to $R_{j}$, in such a way that \eqref{datal22} is
satisfied, i.e.,
\begin{equation}\label{Mj}
   \frac{M_j}{R_{j}}=\la_{0}\equiv \frac{r_{0}}{8\|\chi\|_{L^{2}}}.
\end{equation}
Recalling ~\eqref{gamma1}, we can write for $|s|<c_{0}(N)$
\begin{equation}\label{belowga}
   \left|\int_{0}^{s}
                 \ga''(\sig)d\sig\right|\geq |s|\cdot c_1(N)
\end{equation}
with $c_1(N)>0$. If we choose
\begin{equation}\label{choicemu}
   \mu=c_0(N)/2,
\end{equation}
by \eqref{zeqx}
we obtain
\begin{equation}\label{below5}
   \left|\int_{0}^{\mu z_j(1,x)}
                 \ga''(\sig)d\sig\cdot \chi_{R_j,M_j}(x)\right|\geq
           \frac14c_1(N)c_{0}(N)\chi_{R_j,M_j}(x)
        \quad\hbox{ for }\quad|x|\leq\frac1{R_j}.
\end{equation}
By~\eqref{below5} and~\eqref{below4} we thus get,
recalling \eqref{Mj},
\begin{equation}\label{below6}
   \|\de_t(\uone-\utwo)(1,\cdot)\|^2+\frac{c_0^2}4\|\phi_j\|^2
       \geq\frac18 (c_{0}c_{1})^{2}
       \|\chi_{R_j,M_j}\|^2_{L^2(|x|<1/R_j)}\equiv c_3
\end{equation}
where the constant $c_{3}$ is given by
\begin{equation}\label{c3}
   c_{3}=\frac1{16} (c_{0}c_{1}\la_{0}\|\chi\|_{L^{2}})^{2}
\end{equation}
and is independent of $j$
($\|\chi_{M_{j},R_{j}}\|_{L^{2}}=M_{j}R_{j}^{-1}\|\chi\|_{L^{2}}=
\la_{0} \|\chi\|_{L^{2}} $).

The conclusion of the proof is now quite similar to the general
case $n\geq2$; like before, we choose as data for $v$ the sequences
$$v_0^{(j)}=v_{R_j,M_j,1}(0,x),\qquad
   v_1^{(j)}=\de_t v_{R_j,M_j,1}(0,x)$$
while the data for $w$ are
$$w_0^{(j)}=v_0^{(j)},\qquad
   w_1^{(j)}=v_1^{(j)}+\mu\phi_j    \equiv v_1^{(j)}+c_0\phi_j/2.$$
By \eqref{Mj} the data for $v$ belong to $V$; on the other hand, the
data for $w$ belong to $V$ provided $j$ is large enough, since
$\phi_j\to0$ in $L^2$.
The uniform continuity property implies that
for all $\ep>0$ there exists $\del>0$ such that
$$\|\phi_j\|_{L^2}<\del\qquad\Rightarrow\qquad
    \|\de_t\uone(1,\cdot)-
    \de_t\utwo(1,\cdot)\|_{L^2}\leq\ep$$
and this contradicts ~\eqref{below6}.

\section{Appendix}

The aim of this Appendix is to
prove two multiplicative estimates needed
in the proof of Theorem \ref{thm1}. The first one has the following
form:
\begin{equation}\label{multesta}
      \|fg\|_{H^{s}(\R^{n})}
         \leq C\|f\|_{H^{s}} \cdot\|g\|_{L^{\infty}\cap H^{n/2}}, \ s < n/2.
\end{equation}
Notice that this estimate is \emph{asymmetric} in $f,g.$ We can obtain this estimate
 from the Kato-Ponce  estimate (see Lemma 2.2 in \cite{KPo})
$$\|fg\|_{H^{s}}\leq C\|f\|_{L^{p_1}}\cdot\|J^s g\|_{L^{p_2}}+
        C\|J^s f\|_{L^{p_3}}\cdot\|g\|_{L^{p_4}}$$
which is valid for all $s \geq 0,$ for all
 $p_{2},p_{3}\in]1,\infty[$,  and
$1/2=p_{1}^{-1}+ p_{2}^{-1}= p_{3}^{-1}+ p_{4}^{-1}$;
here $J^{s}=(1-\Delta)^{s/2}$. Then \eqref{multesta} follows
taking $p_{3}=2$, $p_{4}=\infty$,
$$p_{1}=\frac {2n}{n-2s},\qquad
  p_{2}=\frac n{s}$$
and using the Sobolev embeddings
$$\| f\|_{L^{p_{1}}}\leq C\|f\|_{H^{s}},\qquad
  \|J^{s}g\|_{L^{p_{2}}}\leq C\|g\|_{H^{n/2}}.$$

Also the second commutator estimate we need, i.e.,
\begin{equation}\label{below2a}
      \|J^s(fg) - g J^sf \|_{L^2}\leq
         C
      \|f\|_{H^{s}}\cdot\|g\|_{H^{n/2}} ,
      \qquad s < n/2, \qquad n \geq 3,
\end{equation}
can be proved by a similar argument based on
the the Kato-Ponce commutator  estimate (see Lemma 2.2 in \cite{KPo})
$$\|J^s(fg)- g J^s f\|_{L^2}\leq C\|\nabla g \|_{L^{p_1}}\cdot\|J^{s-1} f\|_{L^{p_2}}+
        C\|J^s g\|_{L^{p_3}}\cdot\|f\|_{L^{p_4}}$$
which is valid for all $s \geq 0,$ for all
 $p_{2},p_{3}\in]1,\infty[$,  and
$1/2=p_{1}^{-1}+ p_{2}^{-1}= p_{3}^{-1}+ p_{4}^{-1}$.
Now  \eqref{below2a} follows
taking $p_1=n,$  $p_{3}=s/n$,
$$p_{2}=\frac {2n}{n-2},\qquad
  p_{4}=\frac {2n}{n-2s}$$
and using the Sobolev embeddings
$$\|\nabla g \|_{L^{p_1}}\leq C \|g \|_{H^{n/2}}, \quad \|J^{s-1} f\|_{L^{p_2}}
\leq C \| f\|_{H^s},$$
$$
\|J^s g\|_{L^{p_3}} \leq C \|g \|_{H^{n/2}}, \quad \| f\|_{L^{p_{4}}}\leq C\|f\|_{H^{s}}.$$

For completeness, we  give a self-contained
proof of \eqref{multesta}, \eqref{below2a} and a refined version of \eqref{below2a} involving
homogeneous Sobolev norms; we hope that
our method is of independent interest.

To this end, we must introduce some basic tools from the theory of Sobolev and Besov spaces.

1) \textsc{Difference operators}. Given $h\in\R^{n}$ and
a function $f:\R^{n}\to\C$, we denote
by $f_{j} (x) $ the $j$-th translate of $f$ in the direction $h$:
$$f_{j}(x)=f(x+j\cdot h),\qquad j\in\Z$$
and the \emph{difference operator} $\Del_{h}=\Del$ defined as
$$\Del f=f_{1}-f,\qquad\hbox{i.e.,}\qquad
   \Del f(x)=f(x+h)-f(x) $$
We denote by $\Del^{\ell}$ the iterates of $\Del$.
Trivial properties are $f_{0}\equiv f$, $(f_{i})_{j}=f_{i+j}$,
$\Del^{i}(\Del^{j}f)=\Del^{i+j}f$,
$\Del(f_{j})=(\Del f)_{j}\equiv \Del f_{j}$.

Of special interest here will be the behaviour of the difference operator
with respect to products. We have immediately
$$\Del(fg)=f_{1}g_{1}-fg=f_{1}(g_{1}-g)+(f_{1}-f)g$$
which can be written shortly
$$\Del(f g)=\Del f\cdot g+f_{1}\cdot\Del g.$$
By induction one proves easily the Leibnitz rule
\begin{equation}\label{leibn}
   \Del^{k}(fg)=\sum_{\ell+m=k}
   \genfrac{(}{)}{0pt}{}{k}{\ell}
      \Del^{\ell}f_{m}\Del^{m}g.
\end{equation}

2) \textsc{Sobolev spaces with fractional index}. All the functions
(and the spaces) considered here are defined on the whole
$\R^{n}$. The \emph{homogeneou Sobolev seminorms} $\dot W^{k,p}$ with $k\geq0$ integer, $1<p<\infty$ are defined as
$$\|u\|_{\dot W^{k,p}}
   =\sum_{\al=k}\|D^{\al}u\| _{L^{p}};$$
we write $\dot H^{k}$ for $\dot W^{k,2}$. Thus the standard Sobolev
norms can be written
$$\|u\|_{ W^{k,p}}
   = \|u\|_{L^{p}}+ \|u\|_{\dot W^{k,p}}.$$
For our purposes it is not necessary to enter into the topological details of the definition of the corresponding spaces; only the norms are sufficient, and we shall always apply them to smooth functions.
The $\dot W^{s,p},W^{s,p}$ (semi)norms
with \emph{noninteger} $s>0$ are more troublesome;
the usual definition by interpolation is not well suited to prove
multiplicative estimates. A handier equivalent characterization can
be given using the fractional integrals
$$I_{s,p}(u)=\left(
   \int\!\!\!\int\frac{|\Del_{h}^{[s]+1}u(x)|^{p}}{|h|^{n+sp}}
   dxdh\right)^{1/p}$$
where $[s]$ is the integer part of the noninteger $s>0$, $1<p<\infty$,
and integration is performed over $\R^{2n}$; we shall write
$I_{s,2}=I_{s}$. Then we have
\begin{equation}\label{equiv1}
   \|u\| _{ W^{k,p}}
     \simeq \|u\|_{L^{p}}+ I_{s,p}(u)
\end{equation}
(see e.g. 2.3.1 and  Theorem 2.5.1 in \cite{Triebel78}). The integral $I_{s,p}(u) $ plays
the role of the homogeneous norm; this can be seen by a simple rescaling
argument. For the following application it will be sufficient
to consider the $L^{2}$ case, in which we have a
simple definition using the Fourier transform
$$\|u\|_{\dot H^{s}}=\||\xi|^{s}\widehat u\|_{L^{2}},\qquad
   \|u\|_{H^{s}}\equiv \|(1+|\xi|^{2})^{s/2}\widehat u\|_{L^{2}}
   \simeq \|u\|_{L^{2}}+\|u\|_{\dot H^{s}}.$$
Indeed, let $S_{\la}$ be the scaling operators for $\la>0$
$$(S_{\la}u)(x)=u(\la x);$$
it is easy to check the scaling properties
$$\| S_{\la} u\|_{L^{2}}=\la^{-n/2}\|u\|_{L^{2}},\qquad
     \| S_{\la} u\|_{\dot H^{s}}=\la^{s-n/2}\|u\|_{\dot H^{s}}\qquad
     I_{s}(S_{\la} u)=\la^{2s-n} I_{s}(u). $$
Thus, fixed $u\in  H^{s}$, if we apply the two equivalent
definitions to $S_{\la}u$ we obtain
$$\la^{s-n/2}\|u\|_{\dot H^{s}}+ \la^{-n/2}\|u\|_{L^{2}}\simeq
  \la^{s-n/2} I_{s}(u)+ \la^{-n/2}\|u\|_{L^{2}} .$$
Letting $\la\to\infty$, we obtain immediately
\begin{equation}\label{equiv2}
   \|u\| _{ \dot H^{s}}
     \simeq I_{s}(u) .
\end{equation}

3) \textsc{Besov spaces}. With the same type of norms it is possible
to define the \emph{Besov spaces} $B^{s}_{p,q}$ as follows (see Theorem 2.5.1 in \cite{Triebel78}):
for any $s>0$, $1<p<\infty$, $1<q<\infty$ set
\begin{equation}\label{besov}
   \|f\|_{B^{s}_{p,q}}=\|f\|_{L^{p}}+
     \left(\int
       \frac{\|\Del_{h}^{[s]+1}f\|^{q}_{L^{p}}}{|h|^{n+sq}}dh
     \right)^{1/q}
\end{equation}
and define the spaces accordingly.
From this definition in particular it is evident that
$B^{s}_{p,p}=W^{s,p}$ for noninteger $s$. We shall use the fact that
$$B^{s}_{2,2}\equiv H^{s}$$
for all values of $s$ (including integers).

We finally recall the continuous embedding (see e.g., Theorem 7.58 in \cite{Adams75} and
Theorem 2.8.1 in \cite{Triebel78}): for
$s,t\geq0$
\begin{equation}\label{Sobemb}
      1<p\leq q<\infty,\quad
         s-\frac n p=t-\frac n q\quad\Rightarrow\quad
         W^{s,p}\subseteq W^{t,q}
\end{equation}
and, more generally, the Besov version
\begin{equation}\label{besemb}
   r\in[1,\infty],\quad 1<p\leq q<\infty,\quad
   s-\frac n p=t-\frac n q\quad\Rightarrow\quad
   B^{s}_{p,r}\subseteq B^{t}_{q,r}.
\end{equation}

We are ready to prove our lemma. We use the notation
$$\|u\|_{X\cap Y}= \|u\|_{X}+ \|u\|_{Y}$$
for any two Banach spaces $X,Y$ and $u\in X\cap Y$. We state the following Lemma for smooth functions, the extension to $f,g$ belonging
to the appropriate spaces being obvious.

\begin{lemma}\label{commut}
   For all real $0\leq s<n/2$
   and any smooth functions $f,g$,
   the following inequality holds:
   \begin{equation}\label{multest}
      \|fg\|_{H^{s}(\R^{n})}
         \leq C\|f\|_{H^{s}} \cdot\|g\|_{L^{\infty}\cap H^{n/2}}
   \end{equation}
   and, for all $\la$ with $s<\la<n/2$,
   \begin{equation}\label{multest2}
      \|fg\|_{H^{s}(\R^{n})}
         \leq C\|f\|_{ H^{n/2+s-\la}}\cdot\|g\|_{H^{\la}}.
   \end{equation}
   Moreover, assume that
   \begin{equation}\label{geqc}
      |g(x)|\geq C_{1}>0\hbox{\ \ on the support of $f$;}
   \end{equation}
   then we have also
   \begin{equation}\label{below2}
      \|fg\|_{\dot H^{s}}\geq
         cC_{1}\|f\|_{\dot H^{s}}-c'
      \|f\|_{H^{s}}\cdot\|g\|_{H^{n/2}}
   \end{equation}
   for some constants $c,c'>0$ depending only on $s,n$.
\end{lemma}

\begin{remark}
   Estimate \eqref{multest} can be regarded
   as the limit case of \eqref{multest2} as $\la\to n/2$;
   when $\la\to s$ we obtain
   \eqref{multest} with $f$ and $g$ exchanged.
\end{remark}

\begin{proof}
   Notice that in order to prove \eqref{multest}, \eqref{multest2}
   it is sufficient
   to prove them with the $H^{s}$ norm on the
   left hand side replaced by the
   homogeneous $\dot H^{s}$ norm, since the estimates are trivially true
   for the term $\|fg\|_{L^{2}}$.
   We need two different (but parallel) proofs in the cases $s$ integer
   or noninteger, since we have two different representations of the
   norm in these cases.

   The proof for integer $s$ is is simple. Indeed,
   by the Sobolev embedding
   $$H^{s}(\R^{n})\subseteq L^{\frac{2n}{n-2s}}(\R^{n}),\qquad
         \forall\;0\leq s<\frac n 2$$
   (see \eqref{Sobemb})
   and by H\"older's inequality we have
   \begin{equation}\label{prodl2}
      \|uv\|_{L^{2}}\leq\|u\|_{L^{n/\mu}}\|v\|_{L^{2n/(n-2\mu)}}
        \leq C \|u\|_{H^{n/2-\mu}} \|v\|_{H^{\mu}}
   \end{equation}
   for any real number
   $$0<\mu<\frac n2.$$
   We can apply \eqref{prodl2} to the product of two derivatives
   (here and in the following we shall use the shorthand notation
   $D^{\ell}$ to denote any derivative of order $\ell$):
   \begin{equation}\label{hoso}
      \|D^{\ell}f D^{m}g\|_{L^{2}}\leq
      C\|f\|_{H^{n/2+\ell-\mu}}\|g\|_{H^{m+\mu}}.
   \end{equation}
   For any integer $s<n/2$ we can write
   \begin{equation}\label{hsprod}
     \|fg\|_{\dot H^{s}}\simeq
      \sum_{\ell+m=s}\|D^{\ell}f D^{m}g\|_{L^{2}}.
   \end{equation}
   Now, \eqref{multest2}
   follows directly by applying \eqref{hoso} to each term
   $0\leq\ell<n/2$ with the choice $\mu=\ell+\la-s$, since
   in this case we have $0<\mu<n/2$ for all $\ell=0,...,s$.
   To prove the limit case \eqref{multest}, i.e., with
   $\la=n/2$, the same methods works if we choose
   for $\ell=0,...,s-1$
   $$\mu=\ell+n/2-s$$
   and \eqref{hoso} gives
   $$\|D^{\ell}f D^{m}g\|_{L^{2}} \leq\|f\|_{H^{s}} \|g\|_{H^{n/2}};$$
   but we must consider the term with $\ell=s$ separately
   since $\mu=n/2$ in that case, and we have
   $$\| D^{s} f g\|_{L^{2}}\leq\|g\|_{L^{\infty}}\|f\|_{H^{s}}$$
   and this concludes the proof.

   Consider now
   \eqref{below2} for $s<n/2$ integer; by \eqref{hsprod} we have
   \begin{equation}\label{partdevel}
      \|fg\|_{\dot H^{s}}\geq\|g D^{s}f\|_{L^{2}}-
         c\sum_{\genfrac{}{}{0pt}{}{\ell+m=s}{m\geq1}}
          \| D^{\ell}f D^{m}g\|_{L^{2}}
   \end{equation}
   and applying \eqref{hoso} to each term in the sum, with
   $\mu=n/2-\ell+s$ as above (so that $0<\mu<n/2$), we obtain
   $$\|fg\|_{\dot H^{s}}\geq\|g D^{s}f\|_{L^{2}}-c
      \|g\|_{H^{n/2}} \|f\|_{H^{s}} .$$
   Recalling \eqref{geqc}, we obtain \eqref{below2}.

   From now on, assume $0<s<n/2$ is not an integer. To estimate from
   above $\|fg\|_{\dot H^{s}}$ we use the characterization
   \eqref{equiv2} and the Leibnitz rule \eqref{leibn}:
   $$\|fg\|_{\dot H^{s}}\simeq
     I_{s}(fg)\leq C\sum_{\ell+m=[s]+1}
      \left(\int\!\!\!\int
            \frac{|\Del^{\ell}f_{m} \Del^{m}g |^{2}}{|h|^{n+2s}}
            dxdh\right)^{1/2}. $$
   Thus we need an analogue of \eqref{hoso} for fractional integrals.
   Consider first
   the terms with both $\ell\geq1$ and $m\geq1$.
   By H\"older's inequality we can write for any
   $p^{-1}+q^{-1}=1$ and any $\rho+\sig=s$
   \begin{equation}\label{prodest}
     \int\!\!\!\int
        \frac{|\Del^{\ell}f_{m} \Del^{m}g |^{2}}{|h|^{n+2s}}
        dxdh \leq
     \left(\int\!\!\!\int
        \frac{|\Del^{\ell}f |^{2p}}{|h|^{n+2p\rho}}
        dxdh \right)^{1/p}
    \left(\int\!\!\!\int
        \frac{| \Del^{m}g |^{2q}}{|h|^{n+2q\sig}}
        dxdh \right)^{1/q}
   \end{equation}
   where we replaced $f_{m}$ with $f$ after a translation
   in the variable $x$. The parameters $p,q,\rho,\sig$ must be chosen
   in an appropriate way. First of all we can
   set (since $\ell,m\geq1$)
   \begin{equation}\label{choicerho}
      \rho=\ell-\frac12+\frac{\{s\}}2,\qquad
      \sig=m-\frac12+\frac{\{s\}}2
   \end{equation}
   where $\{s\} =s-[s]$ is the fractional part of $s$, so that
   $$\rho+\sig=s,\qquad \ell=[\rho]+1,\qquad m=[\sig]+1.$$
   Recalling the definition of $I_{s,p}$, from
   \eqref{prodest} we thus obtain
   \begin{equation}\label{prodest2}
     \left(\int\!\!\!\int
        \frac{|\Del^{\ell}f_{m} \Del^{m}g |^{2}}{|h|^{n+2s}}
        dxdh \right)^{1/2}\leq
     I_{\rho,2p}(f)\cdot I_{\sig,2q}(g)\leq
     c\|f\|_{W^{\rho,2p}}\cdot \|g\|_{W^{\sig,2q}}.
   \end{equation}
   Now let $\la$ be such that $s\leq\la\leq n/2$ (extreme
   cases included) and choose $p,q\in]1,\infty[$ as follows
   \begin{equation}\label{choicepq}
      2p=\frac n{\la-\sig},\qquad 2q=\frac{2n}{n-2(\la-\sig)};
   \end{equation}
   notice that $\la-\sig\geq s-\sig>0$ and $n-2(\la-\sig)\geq2\sig>0$.
   Thus by \eqref{Sobemb} we have the embeddings
   \begin{equation}\label{embW}
      H^{n/2+s-\la}\subseteq W^{\rho,2p},\qquad
      H^{\la}\subseteq W^{\sig,2q}.
   \end{equation}
   In conclusion we have
   proved for all $\ell,m\geq1$ with $\ell+m=[s]+1$,
   and any $s\leq\la\leq n/2$, the inequality
   \begin{equation}\label{prodest3}
     \left(\int\!\!\!\int
        \frac{|\Del^{\ell}f_{m} \Del^{m}g |^{2}}{|h|^{n+2s}}
        dxdh \right)^{1/2}\leq
     c\|f\|_{H^{n/2-\la+s}}\cdot \|g\|_{H^{\la}}.
   \end{equation}
   Two terms are left. The term with $m=0$, $\ell=[s]+1$
   is bounded simply by writing
   \begin{equation}\label{prodest4}
     \left(\int\!\!\!\int
        \frac{|g \Del^{[s]+1}f |^{2}}{|h|^{n+2s}}
        dxdh \right)^{1/2}\leq
        \|g\|_{L^{\infty}}\cdot I_{s}(f)\leq
        c \|g\|_{L^{\infty}} \|f\|_{H^{s}}.
   \end{equation}
   On the other hand, the term with $\ell=0$ and $m=[s]+1$ is more
   delicate since we can not use the $L^{\infty}$ norm of $f$. We
   proceed as follows: we apply H\"older
   inequality in $dx$ to obtain
   \begin{equation*}
     \left(\int\!\!\!\int
        \frac{|f_{m} \Del^{[s]+1}g |^{2}}{|h|^{n+2s}}
        dxdh \right)^{1/2}\leq
     \left(\int\!\!\!\int
        \frac{\|f_{m}\|_{L^{n/s}}^{2}
              \|\Del^{[s]+1}g\|_{L^{2n/(n-2s)}}^{2}}{|h|^{n+2s}}
              dh \right)^{1/2}
   \end{equation*}
   and we notice that
   the norm $\|f_{m}\|_{L^{n/s}} = \|f\|_{L^{n/s}} $ is
   independent of $h$ and can be drawn out of the integral. What
   remains is exactly a Besov norm (see \eqref{besov}) and we
   conclude
   \begin{equation}\label{prodest5}
     \left(\int\!\!\!\int
        \frac{|f_{m} \Del^{[s]+1}g |^{2}}{|h|^{n+2s}}
        dxdh \right)^{1/2}\leq
        \| f\|_{L^{2n/(n-2s)}}\cdot\|g\|_{B^{s}_{\frac{n}{s},2}}
        \leq c\|f\|_{H^{s}}\|g\|_{H^{n/2}}
   \end{equation}
   by the  continuous embeddings
   $$H^{s}\subseteq L^{2n/(n-2s)},\qquad
     B^{s}_{\frac{n}{s},2} \subseteq B^{n/2}_{2,2}\equiv H^{n/2}$$
   (see \eqref{Sobemb}, \eqref{besemb}). By \eqref{prodest3} for
   $\la=n/2$, \eqref{prodest4} and \eqref{prodest5} we obtain
   \eqref{multest}.

   By the same method we can write
   \begin{equation*}
     \left(\int\!\!\!\int
        \frac{|f_{m} \Del^{[s]+1}g |^{2}}{|h|^{n+2s}}
        dxdh \right)^{1/2}\leq
     \left(\int\!\!\!\int
        \frac{\|f_{m}\|_{L^{p}}^{2}
              \|\Del^{[s]+1}g\|_{L^{q}}^{2}}{|h|^{n+2s}}
              dh \right)^{1/2}
   \end{equation*}
   where, for an arbitrary $\la$ with $s<\la<n/2$, $p$ and $q$
   are chosen as
   $$2p=\frac n{\la-s},\qquad 2q=\frac{2n}{n-2(\la-s)};$$
   proceeding exactly as in the proof of \eqref{prodest5} we
   obtain
   \begin{equation}\label{prodest6}
     \left(\int\!\!\!\int
        \frac{|f_{m} \Del^{[s]+1}g |^{2}}{|h|^{n+2s}}
        dxdh \right)^{1/2}\leq
        \| f\|_{L^{n/(\la-s)}}\cdot
                \|g\|_{B^{s}_{\frac{2n}{n-2(\la-s)},2}}
        \leq c\|f\|_{H^{n/2+s-\la}}\|g\|_{H^{\la}}.
   \end{equation}
   By \eqref{prodest3} and \eqref{prodest6} we obtain immediately
   \eqref{multest2} for noniteger $s$.

   The proof of \eqref{below2} for noninteger $s$
   proceeds in a similar way. Using again the
   Leibnitz rule \eqref{leibn} we can write
   $$I_{s}(fg)\geq
      \left(\int\!\!\!\int
            \frac{|g\Del^{[s]+1}f|^{2}}{|h|^{n+2s}}
            dxdh\right)^{1/2}
       -c\sum_{\genfrac{}{}{0pt}{}{\ell+m=[s]+1}{m\geq1}}
      \left(\int\!\!\!\int
            \frac{|\Del^{\ell}f_{m} \Del^{m}g |^{2}}{|h|^{n+2s}}
            dxdh\right)^{1/2} $$
   which by \eqref{prodest3} for $\la=n/2$ and
   \eqref{prodest5} implies
   $$I_{s}(fg)\geq
   \left(\int\!\!\!\int
            \frac{|g\Del^{[s]+1}f|^{2}}{|h|^{n+2s}}
            dxdh\right)^{1/2} -
         c \|f\|_{H^{s}}\|g\|_{H^{n/2}}.$$
   Using now assumption \eqref{geqc} we have
   $$I_{s}(fg)\geq C_{1}I_{s}(f)-
       c \|f\|_{H^{s}}\|g\|_{H^{n/2}} $$
   and recalling that $I_{s}(u)\simeq\|u\|_{\dot H^{s}}$,
   we conclude the proof.

\end{proof}

\end{document}